\documentclass[12pt, twoside, leqno]{article}

\usepackage[shortcuts]{extdash}
\usepackage{amsmath,amsthm}
\usepackage{amssymb}

\usepackage{enumerate}

\usepackage{graphicx}

\newtheorem{thm}{Theorem}[section]
\newtheorem{cor}[thm]{Corollary}
\newtheorem{lem}[thm]{Lemma}

\theoremstyle{definition}
\newtheorem{df}[thm]{Definition}

\numberwithin{equation}{section}

\begin{document}


\baselineskip=17pt


\title{On the coincidence of zeroth Milnor-Thurston homology with singular homology}

\author{Janusz Przewocki\\
Institute of Mathematics\\ 
Polish Academy of Sciences\\
00-956 Warszawa, Poland\\
E-mail: j.przewocki@gmail.com
\and 
Andreas Zastrow\\
Institute of Mathematics\\ 
University of Gda\'nsk\\
80-952 Gda\'nsk, Poland\\
E-mail: zastrow@mat.ug.edu.pl}

\date{}

\maketitle


\renewcommand{\thefootnote}{}



\renewcommand{\thefootnote}{\arabic{footnote}}
\setcounter{footnote}{0}


\begin{abstract}
In this paper we prove that the zeroth Milnor-Thurston homology group coincides with singular homology for Peano Continua. Moreover, we show that the canonical homomorphism between these homology theories may not be injective. However, it is proved that it is injective when a space has Borel path-components. 
\end{abstract}
\section{Introduction}
\label{sec:intro}

Milnor-Thurston homology theory was first proposed in Thurston's preprint on \emph{Geometry and Topology of Three-manifolds} \cite[Section 6.1]{Thurst}, with the motivation of providing a more symmetric representation of the fundamental cycle of a  hyperbolic three-manifold than could be achieved with the classical finite chains. It is based on the idea, that when replacing classical finite sums (as they are considered in singular homology theory) by certain measures on the sets of all simplices, then a much more convenient representation of a fundamental
cycle of a hyperbolic manifold can be achieved, in particular, if simplices
with a large hyperbolic volume are to be used. Thus, by its use of  measures
this homology theory may be interpreted in that way, that the classical
finite sums of singular homology theory are replaced by some kind of infinite
sums. However, a canonical homomorphism from singular homology theory to Milnor-Thurston homology theory can always be defined.  

Conversely, the algebraic topology of non-triangulable spaces has produced
some psychological unexpected results, and recently have appeared some
papers \cite{RDiestel1, RDiestel2, RDiestel3}, \cite[Section 1.3]{Georga}  that may be interpreted as a search for a homology-theory
that responds in a more natural way to non-tameness. On the other hand, many of these unexpected results seem to have its origin in the fact, that for defining all classical invariants of algebraic topology only algebraic structures have been used that are based on just considering finite sums and products, while for some spaces like Hawaiian Earrings the topology naturally also allows infinite concatenations. Therefore the idea to investigate how Milnor-Thurston homology theory responds to wildness in topology. 

In \cite{Preprint} the Milnor-Thurston homology groups of the Warsaw Circle were computed, with the surprising result that the zeroth Milnor-Thurston homology group is infinite-dimensional. Milnor-Thurston homology theory satisfies in principle the Eilenberg-Steenrod axioms, but the determination of the isomorphism type of these homology groups (and thus the ``coincidence'' of Milnor-Thurston homology groups with singular homology groups) is only guaranteed for
triangulable spaces. Since the example of the Warsaw Circle (i.e. of a metric compact space) implies that, although zeroth homology is usually related to the number of path components, for non-triangulable spaces the canonical homomorphism from singular to Milnor-Thurston homology can even in this dimension fail to be an isomorphism (in particular: fail to be surjective). There are the following natural two questions:
\begin{itemize}
    \item Is this homomorphism in general injective?
    \item Are there beyond triangulability sufficient criteria, when it will be surjective?
\end{itemize}

In this paper we provide the following answers to these questions:
\begin{itemize}
    \item For Peano-continua we have coincidence, i.e here the canonical homomorphism will always be an isomorphism (cf. Section \ref{sec:peano}).
    \item For spaces with Borel path-components this homomorphism will be at least injective (cf. Section \ref{sec:inject}).
    \item However, we will also provide an example, where it will not even be injective (cf. Section \ref{sec:example}). 
\end{itemize}

Peano Continua are in general not triangulable. Thus
the fact, that the zeroth Milnor-Thurston homology group of a Peano-Continuum will in any case be one-dimensional does neither follow from the Eilenberg-Steenrod Axioms, nor, as the above mentioned example shows,
from the fact that these spaces are path-connected. Nevertheless it holds,
as we will show in this paper.

\section{Preliminaries}
\label{sec:prelim}

We will start this section recalling some facts and definitions concerning Milnor-Thurston homology theory. Then, we will list some results from analysis and measure theory that will be used. 

\subsection{Milnor-Thurston homology theory}

Milnor-Thurston homology theory was defined for differentiable manifolds by Thurston \cite{Thurst} and its initial application was, among others, to prove Gromov's theorem \cite{HJMunkholm}, \cite[Theorem 6.2]{Thurst} that the volume of hyperbolic manifolds is a topological invariant. It was generalised by Hansen 
\cite{SKHansen} to metric spaces and independently by the second author
\cite{AZastrow} to all topological spaces. Its basic definitions can also be found in \cite[\S 11.7]{Ratcl}.

It has been proved that this homology theory satisfies the Eilenberg-Steenrod axioms with a weak version of the Excision Axiom, that is equivalent to the standard version at least for normal spaces \cite[Theorem 4.1]{AZastrow}.  From this fact follows coincidence with singular homology theory for CW-complexes. 

The behaviour of this theory for non-tame spaces is mostly unexplored. Some results in this direction were provided by the second author \cite[Section 6]{AZastrow}, where it is proved that the canonical homomorphism (defined below) between singular homology and Milnor-Thurston homology is not necessarily an isomorphism. Additionally, in \cite{Preprint} it is proved that the first Milnor-Thurston homology group for the Warsaw Circle is trivial, and that the zeroth homology group is uncountable-dimensional, which is an unexpected result.   

Now, we shall briefly present the construction of Milnor-Thurston homology theory. In this paper we use calligraphic letters ($\mathcal{C}$, $\mathcal{H}$, etc.) for constructions in Milnor-Thurston homology theory and ordinary letters for the corresponding constructions in singular homology theory ($C$, $H$, etc.). 

First, we will construct the chain complex $\mathcal{C}_*(X)$ for a given topological space $X$. For that purpose we need to recall some basic notions of measure theory and introduce some notation. 

A $\sigma$-algebra is a family of subsets of some set $\Omega$ such that it is closed under complements and countable unions and it contains the empty set. Any intersection of $\sigma$-algebras is again a $\sigma$-algebra, hence we can consider the smallest $\sigma$-algebra containing a family of subsets $\mathcal{A}$; it shall be denoted $\sigma(\mathcal{A})$ and we say that it is \emph{generated by} $\mathcal{A}$.  

A $\sigma$-additive set function with possibly negative values that is defined on a $\sigma$-algebra is called a \emph{signed measure} if it is zero on the empty set. In this paper we are only interested in the case where $\Omega$ is a topological space and our $\sigma$-algebra is the Borel $\sigma$-algebra (it is the $\sigma$-algebra generated by open sets). Every measure considered here is a signed Borel measure so, for simplicity, we shall call them \emph{measures}. 

A Borel set is called a null-set if all of its Borel subsets have measure zero. A \emph{carrier} of a measure is a set $D$ such that every Borel subset of $\Omega \setminus D$ is a null-set.   Let $\omega \in \Omega$, the measure $\delta_\omega$ is called \emph{the Kronecker measure concentrated on} $\omega$, if its value is one for every set containing $\omega$, and is zero otherwise. Obviously $\{ \omega \}$ is a carrier of such a measure.  



Let $C^0(\Delta^k, X)$ denote the set of singular simplices (continuous functions from the standard simplex $\Delta^k$ to $X$, where $k$ is a non-negative integer). We shall consider $C^0(\Delta^k, X)$ as a topological space equipped  with a compact-open topology. The vector space $\mathcal{C}_k(X)$ of $k$-dimensional chains shall consist of finite measures with a compact carrier (in this paper the notion of compactness does not require Hausdorffness, this is a different terminology than the one used by the second author in \cite{AZastrow}, cf. in particular Definition 1.8 there). 

Given a measurable function $f: \Omega_1 \rightarrow \Omega_2$ and a measure $\mu$ on $\Omega_1$, we can define the \emph{image measure} $f \mu$ by the formula
\begin{displaymath}
(f \mu) (A) = \mu(f^{-1}(A)), \quad \textrm{ for any measurable set A}.
\end{displaymath} 

This construction allows us to define the boundary operator. Indeed, we can see that the natural inclusions of $\Delta^{k-1}$ as a faces of $\Delta^k$ induce  continuous maps on the level of singular simplices. Let $i = 0, 1,...,k$,  we define $\partial_i: \mathcal{C}_{k}(X) \rightarrow \mathcal{C}_{k-1}(X)$ as the image measure construction with respect to the map induced by the inclusion of $i$th face of $\Delta^k$. Finally, the boundary operator is given with the usual formula:
\begin{equation}
\partial = \sum_{i = 0}^k (-1)^i \partial_i. 
\label{eq:differential}
\end{equation}
It was proved \cite[Corollary 2.9]{AZastrow} that $\mathcal{C}_*(X)$ with this boundary operator is a chain complex. 

The Milnor-Thurston homology groups $\mathcal{H}_*(X)$ are then defined as homology groups of this chain complex $\mathcal{C}_*(X)$. Additionally, $\mathcal{C}_*$ can be treated as a functor from the category of topological spaces to the category of chain complexes. Thus, we can define relative homology groups $\mathcal{H}_*(X, A)$ in a natural way.

There is the canonical homomorphism from singular chains to Milnor-Thurston chains (cf. Section \ref{sec:intro})
\begin{eqnarray}
C_k(X; \mathbb{R}) & \rightarrow & \mathcal{C}_k(X), \nonumber \\
\sum_i \alpha_i \sigma_i & \mapsto & \sum_i \alpha_i \delta_{\sigma_i}. \nonumber 
\end{eqnarray}
This homomorphism is a monomorphism if and only if $X$ is $T_0$. Moreover, it maps boundaries to boundaries, thus it induces a homomorphism on the level of homology: 
\begin{displaymath}
H_k(X; \mathbb{R})  \rightarrow  \mathcal{H}_k(X).
\end{displaymath}
It is an isomorphism when $X$ is a CW-complex. 
Additionally, it happens to be an monomorhpism for many wild spaces (e.g. for the zeroth homology of the Warsaw Circle \cite[see proof of Theorem 4]{Preprint} or the second author's example \cite[Section 6]{AZastrow}, the last fact can be proved with the methods of \cite{Preprint}).

\subsection{Results from analysis and measure theory}
In this subsection we define some notions and recall several results that will be used in this paper. 

\begin{thm}
{\normalfont (Hahn \cite[\S 29 Theorem A]{Halmos})} Let $\mu$ be a signed measure on $(\Omega, \mathcal{F})$. Then there exist two disjoint  sets $\Omega^+$, $\Omega^- \in \mathcal {F}$ such that $\Omega = \Omega^+ \cup \Omega^-$ and such that for every $F \in \mathcal{F}$ we have $\mu(F \cap \Omega^+) \geq 0$, $\mu(F \cap \Omega^-) \leq 0$.   
\end{thm}
The decomposition of our space $\Omega$ into sets $\Omega^+$, $\Omega^-$ is not unique. Nevertheless, for two distinct decompositions: $\Omega^+_i$, $\Omega^-_i$, $i = 1, 2$, one can prove that, given any $F \in \mathcal{F}$ we have that $\mu(F\cap \Omega^+_1) = \mu(F\cap \Omega^+_2)$ and $\mu(F\cap \Omega^-_1) = \mu(F\cap \Omega^-_2)$  \cite[p. 122]{Halmos}. Therefore the signed measure $\mu$ can be uniquely decomposed into the following difference of unsigned measures
\begin{displaymath}
\mu = \mu^+ - \mu^-,
\end{displaymath}
where $\mu^+(\cdot) = \mu(\cdot \cap \Omega_+)$, $\mu^-(\cdot) = -\mu(\cdot \cap \Omega_-)$.

\begin{df}
Let $\mu$ be a measure on a space $X$, the variation $|\mu|$  of the measure $\mu$ shall be defined as 
\begin{displaymath}
|\mu|= \mu^+ + \mu^-.
\end{displaymath}
The total variation $\|\mu\|$ shall be defined as
\begin{displaymath}
\|\mu\| = |\mu|(X).
\end{displaymath}
\end{df}

\begin{df}
Let $\mu$ be a signed finite Borel measure. We say that $\mu$ is regular if for every Borel set $B$
\begin{itemize}
    \item $|\mu|(B)$ is the supremum of $|\mu|(K)$ where $K \subset B$ is compact,
    \item $|\mu|(B)$ is the infimum of $|\mu|(U)$ where $U \supset B$ is open.
\end{itemize}  
\end{df}

The space of regular Borel measures on a topological space $X$ shall be denoted $M(X)$. It is a normed space equipped with the total variation norm. Let $C(X)$ denote a space of real continuous functions on a topological space $X$. We have



\begin{thm}
{\normalfont(Compact version of Riesz Representation Theorem \cite[Chapter III, Theorem 5.7]{JBConway})} Let $X$ be a compact Hausdorff space and let $\mu \in M(X)$. Define $F_\mu: C(X) \rightarrow \mathbb{R}$ by:
\begin{displaymath}
F_\mu(f) = \int_{C(X)} f d\mu. 
\end{displaymath}
Then $F_\mu \in C(X)^*$ and the map $\mu \mapsto F_\mu$ is an isometric isomorphism of $M(X)$ onto $C(X)^*$.   
\end{thm}
Here ``$()^*$'' denotes the continuous dual.

We define the following notions as in \cite[p. 41]{PBillingsley}: 
\begin{df}
($\pi$-system) A non-empty family of sets is called a $\pi$-system if it is closed under finite intersections. 
\end{df}
Obviously any topology is a $\pi$-system. 
\begin{df}
($\lambda$-system) A non-empty family of subsets of space $X$ is called $\lambda$-system if it contains $X$, is closed under complements and is closed under countable disjoint unions. 
\end{df}
Notice, that any $\sigma$-algebra is a $\lambda$-system. 

\begin{thm}
{\normalfont(Dynkin's lemma \cite[Theorem 3.2]{PBillingsley})} Let $D$ be a $\lambda$-system and let $P \subset D$ be a $\pi$-system. Then $\sigma (P) \subset D$. 
\label{thm:dynkin}
\end{thm}
\begin{cor}
Let $\mu$ and $\nu$ be Borel measures on a topological space $X$. Suppose $\mu$ and $\nu$ are equal on open sets, then $\mu = \nu$. 
\label{cor:opensets}
\end{cor}
\textbf{Proof.} Let $\mathcal{D}$ be the subset of Borel $\sigma$-algebra such that for every $A \in \mathcal{D}$ we have $\mu(A) = \nu(A)$. We see that $\mathcal{D}$ is a $\lambda$-system. The topology $\tau$ of $X$ is a $\pi$-system such that $\tau \subset \mathcal{D}$. So by Dynkin's lemma we see that $\mathcal{D}$ is in fact the Borel $\sigma$-algebra and hence $\mu = \nu$. 

\begin{flushright}$\square$ \end{flushright}

In construction of measures we shall use the following result of Constantin Carath\'eodory \cite[Theorem 1.3.10]{Prob}:
\begin{thm}
{\normalfont(Carath\'eodory Extension Theorem)} Let $\mu$ be an unsigned measure on an algebra of sets $\mathcal{F}_0$. Then, $\mu$ has a unique extension to a measure on $\sigma(\mathcal{F}_0)$.   
\end{thm}
In fact, if we want to construct a measure it is convenient to define it on some ``smaller'' family of sets:
\begin{df}
We say that a family $\mathcal{S}$ of subsets of $X$ is a semi-algebra if it contains the empty set, it is closed under finite intersections and for any set $E \in \mathcal{S}$ there exists a finite disjoint collection of sets $C_i \in \mathcal{S}$, such that $X \setminus E = \bigcup_i C_i$. 
\end{df}
An example of a semi-algebra over $[-1, 1]$ may be the family of \emph{semi-closed intervals} the of form $[a, b)$ when intersected with $[-1, 1]$.
\begin{cor}
If $\mu$ is a non-negative countably additive set function on a semi-algebra $\mathcal{S}$ such that $\mu(\varnothing) = 0$, then there exists an extension of $\mu$ to $\sigma(\mathcal{S})$. 
\label{cor:carath}
\end{cor}

\textbf{Proof.} The algebra of sets $\mathcal{F}_0$ that is generated by $\mathcal{S}$ has a simple description:
\begin{displaymath}
\mathcal{F}_0 = \left\{ \bigcup_i E_i \mid E_i \textrm{ is a finite collection of subsets of } \mathcal{S} \right\}
\end{displaymath}
It is easy to see that every element of $\mathcal{F}_0$ is in fact a disjoint union of elements in $\mathcal{S}$. Hence $\mu$ has a natural (and well defined!) extension to an additive set function on $\mathcal{F}_0$.

We will prove that it is in fact countably additive. Take a countable collection of subsets $F_j \in \mathcal{F}_0$, such that $F = \bigcup_j F_j \in \mathcal{F}_0$. Each of these sets is a finite disjoint union of elements in $\mathcal{S}$. Namely, $F = \bigcup_i E_i$, $F_j = \bigcup_i E_i^j$. By the intersection property of a semi-algebra we can assume that each $E_i^j$ is a subset of some $E_k$. Thus, we have
\begin{displaymath}
E_i = \bigcup_{E_i^j \subset E_i} E_i^j.  
\end{displaymath} 
Hence, countable additivity   of $\mu$ on $\mathcal{S}$ implies countable additivity of $\mu$ on $\mathcal{F}_0$. Finally, by the Carath\'eodory extension theorem we know that there exists an extension of $\mu$ on $\sigma(\mathcal{F}_0) = \sigma(\mathcal{S})$. 

\begin{flushright} $\square$ \end{flushright}

Let $\mathcal{A}$ and $\mathcal{B}$ be families of subsets of $X$ and let $Y \subset X$, then $Y \cap \mathcal{A}$ denotes $\{Y \cap A \mid A \in \mathcal{A} \}$ and $\mathcal{A} \cup \mathcal{B}$ denotes $\{ A \cup B \mid A \in \mathcal{A}, B \in \mathcal{B} \}$. 
\begin{lem}
Let $A \subset X$ be a subset of a measurable space $(X, \mathcal{F})$. Let $\mathcal{F}$ be generated by a semi-algebra $\mathcal{S}$. Then $A \cap \mathcal{F} = \sigma(A \cap \mathcal{S})$ as a $\sigma$-algebra over $A$.  
\label{lem:trace}
\end{lem}

\textbf{Proof.} The idea of this proof is a slight generalisation of the proof of \cite[Proposition 1.10]{AZastrow} (proofs by this method can also be found in some standard texts on measure theory \cite[I.1 (1.4)]{Bau}, \cite[1.5(Satz 8)]{Hen}). So let $\mathcal{G}$ be the $\sigma$-algebra over $A$ generated by $A \cap \mathcal{S}$. Obviously, we have $\mathcal{G} \subset A \cap \mathcal{F}$. In order to prove the other inclusion, notice that $\mathcal{G} \cup ((X \setminus A) \cap \mathcal{F})$ is a $\sigma$-algebra over $X$ containing $\mathcal{S}$. Thus,  $\mathcal{F} \subset \mathcal{G} \cup ((X \setminus A) \cap \mathcal{F})$. Now, applying to both sides of this inclusion $A \cap $, we obtain $A \cap \mathcal{F} \subset \mathcal{G}$. 

\begin{flushright}$\square$ \end{flushright}  

\begin{lem}
Let $f: X \rightarrow Y$ be a map between a set $X$ and a measurable space $(Y, \mathcal{G})$. Let $\mathcal{G}$ be generated by a semi-algebra $\mathcal{S}$. Then $f^{-1}(\mathcal{F}) = \sigma(f^{-1}(\mathcal{S}))$ as a $\sigma$-algebra over $X$. 
\label{lem:preimage}
\end{lem}

\textbf{Proof.} Without loss of generality we can assume that $f$ is a surjection. This follows from Lemma \ref{lem:trace} and the fact that $f^{-1}(f(X) \cap \mathcal{A}) = f^{-1}(\mathcal{A})$, for every family $\mathcal{A}$ of subsets of $Y$. 

Let $\mathcal{F} \subset f^{-1}(\mathcal{G})$ be the $\sigma$-algebra generated by $f^{-1}(\mathcal{S})$. First, we will prove that $f(\mathcal{F}) := \{f(B) \mid B \in \mathcal{F}\}$ is a $\sigma$-algebra. Countable additivity is proved using good behaviour of images with respect to unions. Finally, let $A = f(B)$ for some $B \in \mathcal{F}$, then $Y \setminus A = f(X \setminus B)$ because $f$ is a surjection and every set in $\mathcal{F}$ is a preimage of a set in $\mathcal{G}$.

We can see that $\mathcal{S} \subset f(\mathcal{F})$, thus $\mathcal{G} \subset f(\mathcal{F})$. Applying the operation $f^{-1}$ to this equation we obtain $f^{-1}(\mathcal{G}) \subset \mathcal{F}$, which proves our lemma.      
\begin{flushright}$\square$ \end{flushright}  

\begin{lem}
Let $G$ be an open set of a metric space $(X, d)$. Then there exists a sequence of continuous functions converging pointwise from below to the characteristic function of $G$.
\label{lem:converg}  
\end{lem} 

\textbf{Proof.} Let $\chi_G$ denote the characteristic function of $G$ and let $f$ be a continuous function on $[0, \infty)$ such that $f(0) = 0$, $f(t) = 1$ for $t \geq 1$ and $0 \leq f \leq 1$. Then $f_n(x) = f (n \cdot d(x, X \setminus G))$ converge pointwise to $\chi_G$ and $f_n \leq \chi_G$ for all $n$. 

\begin{flushright} $\square$ \end{flushright}


\begin{thm}
{\normalfont(Lebesgue Dominated Convergence Theorem \cite[p.229]{HLRoyden})} Let $(X, \mathcal{F}, \mu)$ be a measure space, let $E \in \mathcal{F}$ and let $f_n$ be a sequence of measurable functions on $E$ such that
\begin{displaymath}
|f_n(x)| \leq g(x), \qquad \textrm{for } x \in E
\end{displaymath}
and for an integrable function $g$ on $E$. Suppose
\begin{displaymath}
f_n(x) \rightarrow f(x) 
\end{displaymath}
almost everywhere on $E$. Then,
\begin{displaymath}
\int_E f d\mu = \lim \int_E f_n d\mu.
\end{displaymath}
\label{thm:lebesgue}
\end{thm}



\begin{thm} {\normalfont(Hahn-Banach Theorem \cite[p.187]{HLRoyden})} Let $p$ be a real valued function defined on a vector space $W$ satisfying $p(x+y) \leq p(x) + p(y)$ and $p(\alpha x) = \alpha p(x)$ for all $\alpha \geq 0$. Suppose that $\lambda$ is a linear functional defined on a subspace $V \subset W$ and that $\lambda(v) \leq p(v)$ for all $v \in V$. Then there is a linear functional $\Lambda$ defined on $W$ such that $\Lambda(w) \leq p(w)$ for all $w \in W$ and $\Lambda(v) = \lambda(v)$ for all $v \in V$. 
\end{thm} 

\begin{cor}
Let $W$ be a normed space and let $V \subset W$ be its subspace. Then any bounded linear functional $V \rightarrow \mathbb{R}$ has a bounded extension to $W$.  
\label{cor:hahnbanach}
\end{cor}

\section{Zeroth Milnor-Thurston homology for Peano continua}
\label{sec:peano}

In \cite{Preprint} it has been proved that the Warsaw Circle has uncountable\=/dimensional zeroth Milnor-Thurston homology group. We may suspect that the fact that this space is not locally connected is the reason behind this phenomenon. However, we may notice that there exist path-connected spaces that are not locally path connected and have one-dimensional zeroth homology group. The example may be the Broom Space (it is the cone over the space consisting of the sequence $1/n$ and its limit point).  

Nevertheless, we may ask the opposite question: Does a connected and locally connected space have one-dimensional zeroth Milnor-Thurston homology group? In this section we prove that the answer is affirmative at least when the space is compact (see Theorem \ref{thm:peano}).   

\begin{thm}
Let $f: [0, 1] \rightarrow X$ be a continuous surjection on a metric space $X$. Suppose $\mu$ is a finite Borel measure on $X$, then there exists a measure $\tilde{\mu}$ on $[0, 1]$ such that $f \tilde \mu = \mu$. 
\label{thm:image}
\end{thm}

\textbf{Proof.} Let $V = \{ g \in C([0, 1]) \mid \textrm{there exists } h \in C(X) \textrm{ such that } g = h \circ f \}$. We see that $V$ is a nonempty linear space. Let $g \in V$, thanks to surjectivity of $f$ the function $h \in C(X)$ such that $g = h \circ f$ is unique. We shall denote it by $h_g$. Notice, that $h_g$ is linear with respect to $g$. 

One can show that the linear functional below is bounded (it follows from the fact that the norm on $V$ is supremum norm and that $\mu$ is finite)
\begin{displaymath}
g \mapsto \int_X h_g d \mu. 
\end{displaymath}
By Corollary \ref{cor:hahnbanach} there exists a bounded extension $\xi$ of this linear functional. Then, by Riesz Representation Theorem we know that there exists a Borel measure $\tilde \mu$ such that
\begin{displaymath}
\xi(g) = \int_{[0, 1]} g d \tilde \mu. 
\end{displaymath} 

Now, we shall prove that $f \tilde \mu = \mu$. By Corollary \ref{cor:opensets} it suffices to check this only for an arbitrary open set $G \subset X$. By Lemma \ref{lem:converg} there exists a sequence $(h_n)_{n \in \mathbb{N}}$ of positive functions that is pointwise convergent to $\chi_G$ and such that $h_n \leq \chi_G$. Let $g_n = h_n \circ f$. Then for each $n$ the function $g_n \in V$, and the sequence $(g_n)_{n \in \mathbb{N}}$ is pointwise convergent from below to $\chi_{f^{-1}(G)}$. 

We know that 
\begin{displaymath}
\int_{[0, 1]} g_n d \tilde \mu = \xi(g_n) = \int_X h_n d \mu .
\end{displaymath}
Using Theorem \ref{thm:lebesgue} on the both sides of the above equation we get
\begin{displaymath}
\int_{[0, 1]} \chi_{f^{-1}(G)} d \tilde \mu  = \int_X \chi_G d \mu, 
\end{displaymath}
which means that $\tilde \mu ( f^{-1}(G)) = \mu(G)$, hence $f \tilde \mu (G) = \mu(G)$. 

\begin{flushright} $\square$ \end{flushright} 

\textbf{Remark.} By convention the $i$th face of a simplex is the face opposed to the $i$th vertex. Consequently, with respect to the boundary sides as
defined in (\ref{eq:differential}) for a 1-simplex $\sigma : [0,1] \rightarrow X$, we obtain that $\partial_0(\sigma)$  maps to   $\sigma(1)$ and $\partial_1(\sigma)$  maps to   $\sigma(0)$.

\begin{thm} 
If $X$ is a Peano continuum, then $\mathcal{H}_0(X) \cong \mathbb{R}$.    
\label{thm:peano}
\end{thm}

\textbf{Proof.} A Peano continuum is a metric continuum that is locally connected. We shall use the Hahn-Mazurkiewicz theorem \cite[Theorem 3-30]{JGHocking} which states that there is a continuous surjection $f: [0, 1] \rightarrow X$. Let $\mu \in \mathcal{C}_0(X)$ represent some homology class. From Theorem \ref{thm:image} we know that there exists a measure $\tilde \mu$ on $[0, 1]$ such that $f \tilde \mu = \mu$. 

Next, let us define $g: [0, 1] \rightarrow C^0(\Delta^1, X)$ with the following formula: $g(x)(t) = f(t x)$. Let $\nu = g \tilde \mu$, we shall prove that $\partial \nu = \mu - \mu(X) \delta_{f(0)}$. Take any Borel subset $A \subset X$, then  
\begin{equation}     
\partial \nu (A) = \tilde \mu(g^{-1} (\partial_0^{-1}  A)) - \tilde \mu(g^{-1}
(\partial_1^{-1} A)).
\label{eq:boundary}
\end{equation}
Suppose $f(0) \notin A$. Then, $g^{-1} (\partial_0^{-1}  A) = f^{-1}(A)$   and
$g^{-1} (\partial_1^{-1} A)$  is empty, so equation (\ref{eq:boundary}) reduces to:
\begin{displaymath}
\partial \nu (A) = \tilde \mu (f^{-1}(A)) = \mu(A).
\end{displaymath}
And when $f(0) \in A$, we have $g^{-1} ( \partial_0^{-1}  A) = f^{-1}(A)$   and
 $g^{-1} (\partial_1^{-1} A) =  [0, 1]$, then equation (\ref{eq:boundary}) reduces to:
\begin{displaymath}
\partial \nu (A) = \tilde \mu (f^{-1} (A)) - \tilde \mu(f^{-1} (X)) = \mu(A) - \mu(X).
\end{displaymath}
From that, we see that every cycle $\mu \in \mathcal{C}_0(X)$ is homological to the measure $\mu(X) \delta_{f(0)}$. 

The Kronecker measure $\delta_{f(0)}$ is non-trivial on the level of homology. Indeed, to the contrary suppose $\partial \alpha = \delta_{f(0)}$ for some measure $\alpha$. By the obvious fact that every singular 1-simplex in $X$ has both its endpoints in $X$ we have the following equality between sets: $\partial_0^{-1} X$ = $\partial_1^{-1} X$. Hence, $(\partial \alpha) (X) = \alpha(\partial_0^{-1} X) - \alpha(\partial_1^{-1} X) = \alpha(\partial_1^{-1} X) - \alpha(\partial_1^{-1} X) = 0$. That contradicts the fact that $\delta_{f(0)}(X) = 1$. Thus, our zeroth homology group is a one-dimensional vector space. 


\begin{flushright} $\square$ \end{flushright}

\section[Is the canonical map a monomorphism?]{Is the canonical map from singular homology to Milnor-Thurston homology a monomorphism?}
\label{sec:inject}

In Section \ref{sec:prelim} we have seen that there exists a canonical homomorphism from singular homology groups to Milnor-Thurston homology groups
\begin{displaymath}
H_k(X; \mathbb{R}) \rightarrow \mathcal{H}_k(X),
\end{displaymath}
where $X$ is a topological space and $k$ is a non-negative integer. 

The second author \cite{AZastrow} showed that this canonical homomorphism is an isomorphism when $X$ is a CW-complex. So in this case it is obviously an injection. Moreover, in \cite{Preprint} it was observed that it is an injection, when $X$ is the Warsaw Circle. 

In this section we consider the question whether this homomorphism is always an injection for $k = 0$. We shall prove the following theorem
\begin{thm} 
Let $X$ be a topological space with Borel path-components. Then, the canonical map $H_0(X; \mathbb{R}) \rightarrow \mathcal{H}_0(X)$ is an injection.    
\label{thm:inject}
\end{thm}
\begin{lem} 
Let $X$ be a topological space with Borel path-components. Let $\mu$ be a measure on $C^0(\Delta^1, X)$, such that $\partial \mu = \nu_{X_1} - \delta_{x_0}$, where $\nu_{X_1}$ is concentrated on a set $X_1 \subset X$ and $x_0 \notin X_1$. Then there exists a path starting at $x_0$ with its endpoint in $X_1$.   
\end{lem}

\textbf{Proof.} Let $Y$ be the path-component containing $x_0$. Notice that $\partial_0^{-1}(Y) = \partial_1^{-1}(Y)$. Thus, we have
\begin{displaymath}
(\partial \mu) (Y) = \mu (\partial_0^{-1}(Y)) - \mu (\partial_1^{-1}(Y)) = 0. 
\end{displaymath}
Now, assume that there is no path from $x_0$ to any point of $X_1$. That is, $X_1$ intersects $Y$ in the empty set. As a consequence, $(\partial \mu)(Y) = -1$ which contradicts the above calculations.

\begin{flushright} $\square$\end{flushright}

\textbf{Proof of Theorem \ref{thm:inject}.}  Suppose, we have a singular cycle $z = \sum_{i = 1}^k \alpha_i x_i$ such that $z = \partial \mu$ for some $\mu \in \mathcal{C}_1(X)$. We will construct a singular chain with that property. We will proceed inductively. 
So assume this fact is true for cycles with number of simplices less than $k$. We will use the above lemma to construct a path from $x_k$ to some point of $X_1 = \bigcup_{i = 1}^{k-1} \{x_i\}$. The measure $\mu / \alpha_k$ satisfies the assumptions of the above lemma, so there exists a path $\sigma_k$  connecting $x_k$ to, say, $x_j$. Let $\tilde \mu = \mu + \alpha_k \delta_{\sigma_k}$, we can see that $\partial \tilde \mu$ has $k-1$ simplices. So, there is a singular chain $c$ such that $\partial c = \partial \tilde \mu$. Now, we see that $\partial (c - \alpha_k \delta_{\sigma_k}) = \sum_{i = 1}^k \alpha_i x_i$ which ends our proof. 
\begin{flushright} $\square$\end{flushright}

\section{A space with a non-injective canonical homomorphism}
\label{sec:example}

The assumption that $X$ has Borel path components was crucial in the proof of the previous theorem. Now, we will construct a counterexample showing that this assumption cannot be omitted. Namely, we will construct a topological space $X$, where there exists a measure $\nu \in \mathcal{C}_1(X)$ such that $\partial \nu = \delta_{x_1} - \delta_{x_0}$ where the points $x_1, x_0 \in X$ lie in separate path components.  

The following lemma will allow us to perform our construction:
\begin{lem}
There exists a partition $[-1, 1] = A \cup B$, where $A$ and $B$ are not Lebesgue measurable and every Borel subset of $A$ or $B$ is of measure zero.  
\end{lem}

\textbf{Remark.} In principle the existence of such sets follows from the ergodicity of the action of rational numbers on the reals \cite[Chapter I, Prop. 4.5.1]{GMargulis}. However, we came to the conclusion that the statement that we need in our context is easier accessible by the following elementary proof than by a quote to a big theory.

\textbf{Proof.} First, we will find such a partition for $S^1 = \mathbb{R} / \mathbb{Z}$. It is enough to show that there exists a set $A \subset S^1$ with Lebesgue inner measure zero and full Lebesgue outer measure (here we normalise the Lebesgue measure $\lambda$ in a way that $\lambda(S^1) = 2$). Indeed, if we have $\lambda_*(A) = 0$ and  $\lambda^*(A) = 2$, then the set $B$ can be defined as a complement of $A$. We see that
\begin{displaymath}
\lambda_*(B) = \sup_{B \supset O \in \mathcal{B}(S^1)} \lambda(O) = \sup_{A \subset O' \in \mathcal{B}(S^1)} (2 - \lambda(O')) = 2 - \lambda^*(A) = 0, 
\end{displaymath}
thus every Borel subset of $B$ has indeed Lebesgue measure zero.

In order to construct the subset $A$, we will use the natural action of  $G = \mathbb{Q} / \mathbb{Z}$ on $S^1$. It is known that $\mathcal{B}(S^1)$ has the cardinality of continuum \cite[Theorem 3.3.18]{SrivastavaSM}. Let $(B_\alpha)_{\alpha < \mathfrak{c}}$ denote the family $\mathcal{B}(S^1)$ with a well-ordering. This well-ordering exists by the well-ordering theorem, which is equivalent to the Axiom of Choice. Using transfinite induction, we shall construct a sequence of elements $(x_\alpha)_{\alpha < \mathfrak{c}}$. 

Suppose, we have chosen $x_\beta$ for all $\beta < \alpha$. Then, we chose $x_\alpha$ that satisfy the following conditions:
\begin{itemize}
    \item for every $\beta < \alpha$, the element $x_\alpha$ lies in a different orbit of $G$-action than $x_\beta$,
    \item if complement of $B_\alpha$ is uncountable, then $x_\alpha \in S^1 \setminus B_\alpha$. 
\end{itemize} 
Elements satisfying both of these conditions always exist. That is because, the number of $G$-orbits is continuum. Moreover, if $\kappa$ denote the number of $G$-orbits that intersect $S^1 \setminus B_\alpha$ in a nonempty set, then the cardinality of $S^1 \setminus B_\alpha$ is less then $\aleph_0 \cdot  \kappa = \max(\aleph_0, \kappa$). Thus, if cardinality of $S^1 \setminus B_\alpha$ is uncountable then it is continuum, which is true for every uncountable Borel set \cite[Theorem 3.2.7]{SrivastavaSM}. Consequently, we see that $\kappa = \mathfrak{c}$, so there are continuum-many orbits we can choose the element $x_\alpha$ from.  

Now, we shall prove that the set $A := \{x_\alpha\}_{\alpha < \mathfrak{c}}$ has the desired properties. Suppose, we have a Borel set $O \subset A$, then both $A$ and $O$ intersect each orbit of $G$ in a set with at most one element. From that, the family $G + O := \{g+O \mid g \in G\}$ consist of pairwise disjoint sets. Now, suppose $\lambda(O) > 0$, then $\lambda\left(\bigcup (G+O)\right) = \sum_{g \in G} \lambda(g+O) = \infty$, which is impossible. Hence, $\lambda(O) = 0$. On the other hand, consider $O \supset A$. If $O$ has a countable complement, then it has full Lebesgue measure. Otherwise, from the fact that $O = B_\alpha$ for some $\alpha < \mathfrak{c}$, we know that $x_\alpha \notin O$, which contradicts $O \supset A$. 

Finally, we can construct our decomposition of the interval $[-1, 1]$. There exists a continuous, measure preserving, map $f: [-1, 1] \rightarrow S^1$ which identifies both ends of the interval. In order to get a partition of $[-1, 1]$ we take preimages of $S^1 = A \cup B$. The properties of the partition are conserved, since $f$ preserves measure. 





\begin{flushright} $\square$ \end{flushright} 

Now, we will start our construction. Take $N = [-1, 1] \setminus \mathbb{Q}$ with the topology induced from the real line. By the above theorem, there exist disjoint non-Lebesgue measurable sets such that $N = N_0 \cup N_1$ and for any Borel set $A \subset N_i$ we have $\lambda(A) = 0$. 

In order to get two connected components, the next stage of our construction will be taking cones over $N_0$ and $N_1$. So, identify $N$ with the subset of $\mathbb{R} \times \{0\} \subset \mathbb{R}^2$. We define the cone $C N_0$ as the union of affine intervals connecting the points of $N_0$ with $x_0 := (0, 1)$. Analogously, let $C N_1$ be the union of intervals connecting $N_1$ with $x_1 := (0, -1)$. 

Notice, that the above construction of a cone is different than usual. Taking the Cartesian product with the interval, and then collapsing one face to a point yields a different neighbourhood system of the cone-point than the one induced from the plane. 



Let $Y := C N_0 \cup C N_1$ and let $I_0$, $I_1$ be disjoint copies of $[0, 1]$ and $[-1, 0]$, respectively. Let us identify the point $1$ of $I_0$ with the vertex $x_0 \in Y$ and point $-1$ of $I_1$ with $x_1 \in Y$. This is the underlying  set of our space $X$.

The topology on $Y$ is induced from $\mathbb{R}^2$. Thus, by choosing a neighbourhood basis of each point of $I_i$ for $i = 0, 1$, we will complete the definition of the topology of $X$. 

Let $\mathcal{J}_i$ denote the family of finite subsets of $N_i$. Then for each $J \in \mathcal{J}_i$ let $C N_i ^ J$ denote the sub-cone $ C(N_i \setminus J) \subset CN_i$.  


Now, let $t \in I_i$ (remember that $t$ is identified with a real number). Its basis of neighbourhoods shall be
\begin{displaymath}
\mathcal{B}_t = \{U^\varepsilon \cup U_{J, K}^\varepsilon \mid \varepsilon > 0, J \in \mathcal{J}_i, K \in \mathcal{J}_j\}     
\end{displaymath}
where $U^\varepsilon = (t-\varepsilon, t+\varepsilon) \cap I_i$ and $U_{J, K}^\varepsilon = \{(x, y) \in \mathbb{R}^2 \mid t - \varepsilon < y < t + \varepsilon \} \cap (C N_i^J \cup C N_j^K)$, for $j = 1-i$.  

Let $y_i$ denote the endpoint of $I_i$ different form $x_i$ and let $T = N \cup \{y_0, y_1\}$ with the topology induced from $X$. 
\begin{lem}
Every continuous map $f: [0, 1] \rightarrow T$ is constant.
\end{lem}

\textbf{Proof.} The lemma is true if $f([0, 1]) \subset N$. So, suppose that $f^{-1}(\{y_0, y_1\})$ is not empty.  


First consider the case when $f^{-1}(N)$ is empty. Then $[0, 1]$ can be decomposed into the disjoint union of closed sets: $f^{-1}(\{y_0\}) \cup f^{-1}(\{y_1\})$, this contradicts connectivity of $[0, 1]$.  

Next, let $f^{-1}(N)$ be nonempty. Notice that it is an open set because $N$ is open in $T$. Therefore, it must be a countable disjoint union of open nonempty intervals. Now, take $(a, b)$ to be one of these intervals. By assumption, $f(a) = y_i$ for some $i$. Because $(a, b)$ is connected $f$ should be constant on it with a value, say, $x \in N$. There exists a neighbourhood of $y_i$ without $x$, therefore $f$ is discontinuous at $a$. 
\begin{flushright} $\square$ \end{flushright}

\begin{lem}
The points $x_0$ and $x_1$ lie in different path-components. 
\label{lem:pathcomp}
\end{lem}

\textbf{Proof.} 
Suppose that there is a path $\alpha: [0, 1] \rightarrow X$ connecting $x_0$ and $x_1$. 
Notice that there is a supremum $t_0$ of points $t$ such that $\alpha(t) = x_0$. From the continuity of $\alpha$ we see $\alpha(t_0) = x_0$. Similarly, there exists a infimum $t_1$ of points $t > t_0$ such that $\alpha(t) = x_1$. Now, we have that the points between $t_0$ and $t_1$ are mapped into $X \setminus \{x_0, x_1\}$.
 
Take a point $a \in [t_0, t_1]$ close enough to $t_0$ so that $\alpha(a) \in C N_0$ and take a point $b \in [t_0, t_1]$ close enough to $t_1$ so that $\alpha(b) \in C N_1$. We see that the interval $[a, b]$ is mapped into $X \setminus \{x_0, x_1\}$, so we can construct a path $\beta: [0, 1] \rightarrow X \setminus \{x_0, x_1\}$ connecting a point of $C N_0$ with a point of $C N_1$. 

There is the obvious retraction $r: X \setminus \{x_0, x_1 \} \rightarrow T$ that maps each point to the end-point of its ray in the respective cone. By the above lemma the function $r \circ \beta$ is constant, hence $\beta$ maps the interval $[0, 1]$ into a single ray of one of the cones. Consequently, it cannot connect points in separate cones. 

\begin{flushright}$\square$ \end{flushright} 




Now, we shall construct our measure $\nu$ on $C^0(\Delta^1, X)$. It will consist of two parts, one concentrated on simplices in $C N_0$ and the other concentrated on simplices in $C N_1$. Their carriers shall consist of simplices connecting points of $N$ with the respective vertex.  

To get a convenient description of carriers for our measures we shall still treat $Y$ as a subset of $\mathbb{R}^2$ (in the way described above). Let $\sigma_{x_0}^x$ be the singular simplex such that $\sigma_{x_0}^x (t) = (f_x(t), 1-t)$, where $f_x$ is the unique affine function such that $\sigma_{x_0}^x(t) \in Y$ and $\sigma_{x_0}^x (1) = x$. In the analogous way we define simplex $\sigma_x^{x_1}$ for $x \in N_1$ (the direction is such that $\sigma_x^{x_1}(0) = x$). 

Now, our carriers shall be $S_0 = \{\sigma_{x_0}^x \in C^0(\Delta^1, X) \mid x \in N_0\}$ and $S_1 = \{\sigma_x^{x_1} \in C^0(\Delta^1, X) \mid x \in N_1\}$. 

Notice that each of $S_i$ is not compact, however if we add to $S_i$ the respective paths connecting $x_i$ with $y_i$ (parametrised in the proper way) we shall get a compact sets of simplices. Thus, our measure $\nu$ shall have a compact carrier.  

\begin{lem}
$S_0$ and $S_1$ are Borel sets in $C^0(\Delta^1, X)$. 
\label{lem:borel}
\end{lem}


\textbf{Proof.} First, we will show that it is sufficient to prove that $S_i$ are Borel in $C^0(\Delta^1, Y)$. To do this we show that $C^0(\Delta^1, Y)$ is Borel in $C^0(\Delta^1,X)$. The desired conclusion follows, since every Borel subset of a Borel subspace is Borel in the bigger space. 

Take $i = 0, 1$, and let $U_n^i$ denote a sequence of neighbourhoods of $x_i$ such that $\bigcap_n U_n^i = \{x_i\}$. Now, let  $Y_n = Y \cup U_n^0 \cup U_n^1$. We see that each $Y_n$ is an open set in $X$ and $\bigcap_n Y_n = Y$. By this fact and the definition of the compact-open topology, $C^0(\Delta^1, Y_n)$ is open in $C^0(\Delta^1, X)$. The intersection of $C^0(\Delta^1, Y_n)$ is $C^0(\Delta^1, Y)$, so it is a Borel set.  

Now, we shall prove that the each of $S_i$ is closed in $C^0(\Delta^1, Y)$. The space $C^0(\Delta^1, Y)$ is metrizable, thus it is enough to show that both $S_i$ contain limit points of all sequences. Let $\sigma_n$ be a sequence of singular 1-simplices in $\mathbb{R}^2$ with affine parametrisation, say, $\sigma_n(t) = (a_n + b_n t, c_n + d_n t)$. Such a sequence is convergent iff sequences of coefficients $a_n, b_n,$ etc. are convergent. 

Now, take a sequence of 1-simplices $(\sigma_n) \subset S_0 \subset C^0(\Delta^1, Y) \subset C^0(\Delta^1, \mathbb{R}^2)$ convergent in  $C^0(\Delta^1, Y)$. By the above observation a limit of such a sequence is a 1-simplex with affine parametrisation that connects $x_0$ with a point of $N$. However, any such simplex is an element of $S_0$, hence $S_0$ is closed. Analogously, we prove that $S_1$ closed. 

\begin{flushright}$\square$\end{flushright}

We preferred to state the following lemma in an abstract way. Its assumptions are satisfied in our case. Namely, take $Z = C^0(\Delta^1, Y)$,  $f_i = \partial_i$, $M_i = S_i$ (this yields $R_i = N_i$), and $S_i$ is homeomorphic to $N_i$ (cf. Lemma \ref{lem:homeo1} and Lemma \ref{lem:homeo2}).  

\begin{lem}
Let $f_i: Z \rightarrow [-1, 1]$ for $i = 0, 1$ be continuous functions on a topological space $Z$ such that there exist disjoint Borel subsets $M_i$ of $Z$ with the following properties:
\begin{itemize}
    \item Every $R_i = f_i(M_i)$ is dense in $[-1, 1]$,
    \item $R_i$ are disjoint
    \item $R_0 \cup R_1$ is a full-measure Borel set 
    \item Every Borel subset of $f_i(R_i)$ has Lebesgue measure zero,
    \item $f_i|{M_i}$ is a homeomorphism of $M_i$ and $R_i$. 
\end{itemize}
Then
\begin{enumerate}
    \item Every Borel set in $M_i$ has the form $M_i \cap f_i^{-1}(B)$ for some Borel subset of $[-1, 1]$, 
    \item The semi-algebra $\mathcal{I}_i = \{M_i \cap f_i^{-1}(I) \mid I \subset [-1, 1] \textrm{ is a semi-closed in-}$ $\textrm{terval } \} $ generates Borel subsets of $M_i$,  
    \item The set functions $\mu_i: M_i \cap f_i^{-1}(I) \mapsto \lambda(I)$, where $\lambda$ denotes the Lebesgue measure and $I$ is a semi-closed subinterval of $[-1, 1]$, can be extended to a Borel measure $\mu_i$ on $M_i$.  
\end{enumerate}
\label{lem:mainlemma}
\end{lem}

\textbf{Proof.} To prove the first statement take a Borel subset $A$ of $M_i$. Then $f_i(A)$ is a Borel subset of $R_i$. Notice, that every Borel subset of $R_i$ is an intersection of Borel subset of $[-1, 1]$ and $R_i$, which proves the first statement. 

To prove the second statement we need to notice that  $\mathcal{I}_i = \{M_i \cap f_i^{-1}(I) \mid I \subset [-1, 1] \textrm{ is a semi-closed interval} \} $ is a semi-algebra. Then Lemma \ref{lem:preimage} and Lemma \ref{lem:trace} combined together give us this result.  
 
In order to prove the third statement it is enough to show that $\mu_i$ are countably additive (see Corollary \ref{cor:carath}). So, let us take a pairwise disjoint countable family $\{ A_j = M_i \cap f_i^{-1}(I_j)  \} \subset M_i \cap f_i^{-1}(\mathcal{I})$, such that the union of this family is  $A \in M_i \cap f_i^{-1}(\mathcal{I})$. Thus, the set $A$ is of the form $M_i \cap f_i^{-1}(I)$ for some semi-closed interval $I$. 

We claim that $\{I_j\}$ is a pairwise disjoint family. To the contrary, assume that two of these sets, say, $I_1$ and $I_2$, have the non-empty intersection $[a, b)$, for some real numbers $a < b$. Consequently, $A_1 \cap A_2 = M_i \cap f_i^{-1}([a, b))$ and is non-empty, since $R_i$ is dense in $[-1, 1]$. However, our family of sets is disjoint, hence we got a contradiction. 


Moreover, we claim that $I \setminus \bigcup_j I_j$ is a Borel subset of $[-1, 1] \setminus R_i$. Indeed, from the fact that $A$ is the union of $A_j$, we get $M_i \cap f_i^{-1}(\bigcup_j I_j) = M_i \cap f_i^{-1}(I)$. 

Next, we see that $[-1, 1] \setminus R_i = R_k \cup ([-1, 1] \setminus R_0 \cup R_1 )$ for $k = 1 - i$. Thus, $I \setminus \bigcup_j I_j$ can be decomposed into two parts. The first one is a Borel subset of $[-1, 1] \setminus R_0 \cup R_1$ and hence it is a null-set as a subset of a null-set. The second one is a Borel subset of $R_k$ and every such subset is a null-set. As a consequence $I \setminus \bigcup_j I_j$ is a null-set, which yields $\lambda(I) = \sum_j \lambda(I_j)$. This fact proves that $\mu_i$ is countably additive. 


\begin{flushright}$\square$\end{flushright}

\begin{lem}
$\partial_0|_{S_0}: S_0 \rightarrow N_0$ is a homeomorphism.     
\label{lem:homeo1}
\end{lem}

\textbf{Proof.} 
It is obviously a bijection; in order to prove that its inverse is continuous, take a convergent sequence of $x_n \in N_0$ with a limit $x \in N_0$. Observe that for each $n$ the simplex $(\partial_0|_{S_0})^{-1} (x_n)$ connects $x_0$ with $x_n$. Because all simplices have affine parametrisation, the limit of $(\partial_0|_{S_0})^{-1} (x_n)$ is the simplex connecting $x_0$ with $x$ that has affine parametrisation. Hence $(\partial_0|_{S_0})^{-1}$ is continuous.

\begin{flushright} $\square$ \end{flushright}

\begin{lem}
$\partial_1|_{S_1}: S_1 \rightarrow N_1$ is a homeomorphism.
\label{lem:homeo2}     
\end{lem}

\textbf{Proof.} Is analogous as in the previous case. 

\begin{flushright} $\square$ \end{flushright}

Now, let $\nu_i$ be the measure on Borel subsets of $S_i$ that exists by Assertion 3 of Lemma \ref{lem:mainlemma}. We can extend the  measures $\nu_i$ for $i = 0, 1$ to the Borel $\sigma$-algebra of $C^0(\Delta^1, X)$ with the formula
\begin{displaymath}
\nu_i(A) = \nu_i(A \cap S_i), \quad \textrm{ for any Borel subset } A \textrm{ of } C^0(\Delta^1, X), 
\end{displaymath}
which is well-defined thanks to Lemma \ref{lem:borel}. 

Now, let us define $\nu = \nu_1 + \nu_0$. Finally, we can prove our result.
\begin{thm}
The canonical homomorhpism $h: H_0(X; \mathbb{R}) \rightarrow \mathcal{H}_0(X)$ is not a monomorphism. 
\end{thm}

\textbf{Proof.} The singular homology class of the cycle $z = x_1 - x_0$ is nontrivial in $H_0(X; \mathbb{R})$ since $x_0$ and $x_1$ lie in separate path components (see Lemma \ref{lem:pathcomp}). The canonical homomorphism maps this class to the class of the cycle $\delta_{x_1} - \delta_{x_0}$ in $\mathcal{H}_0(X)$, where $\delta$ denotes the Kronecker measure. We shall prove that the homology class of this measure is trivial. In fact, we shall show that for the measure $\nu$ defined above we have
\begin{displaymath}
\partial \nu = 2(\delta_{x_1} - \delta_{x_0}).
\end{displaymath}

The crucial step of our proof is to show that every Borel subset of $N$ is of $\partial \nu$-measure zero. So, let $B \subset N \subset [-1, 1]$ be a Borel set. Notice, that $\nu_1(\partial_0^{-1}(B)) = 0$ because $S_1 \cap \partial_1^{-1}(B)$ is empty. Similarly, $\nu_0(\partial_1^{-1}(B)) = 0$. As a consequence we see  
\begin{displaymath}
(\partial \nu)(B) = \nu_0(\partial_0^{-1}(B)) - \nu_1(\partial_1^{-1}(B)).
\end{displaymath} 

Now, notice that if $B = I \cap N$, where $I$ is an interval, we have $(\partial \nu)(B) =  \nu_0(\partial_0^{-1}(I)) - \nu_1(\partial_1^{-1}(I)) = \lambda(I) - \lambda(I) = 0$. So the $\lambda$-system of Borel sets that satisfy $(\partial \nu)(B) = 0$ contains a semi-algebra generating Borel subsets of $N$. Every semi-algebra is a $\pi$-system, so by Theorem \ref{thm:dynkin} we see that our assertion is true for every Borel set. 

\begin{flushright} $\square$ \end{flushright}


\begin{thebibliography}{999}

\bibitem{Prob} R. B. Ash, \emph{Probability and measure theory}, San Diego [etc.]: Academic Press, 2000

\bibitem{Bau} H. Bauer, \emph{Wahrscheindlichkeitstheorie und Grundz\"uge der Ma\ss theorie}, 2 Auflage, Berlin-New York: deGruyer,  1974 

\bibitem{PBillingsley} P. Billingsley, \emph{Probability and measure}, 3rd ed., New York [etc.]: John Wiley \& Sons, 1995

\bibitem{JBConway} J. B. Conway, \emph{A course in functional analysis}, 2nd ed, New York [etc.]: Springer-Verlag, 1990

\bibitem{RDiestel1} R. Diestel, P. Spr\"ussel, \emph{Locally finite graphs with ends: a topological approach, III, Fundamental group and homology},  Discrete
Math. 312 (2012), no. 1, 21--29

\bibitem{RDiestel2} R. Diestel, P. Spr\"ussel, \emph{On the homology of locally
compact spaces with ends}, Topology Appl. 158 (2011), no. 13, 1626--1639


\bibitem{RDiestel3} R. Diestel, P. Spr\"ussel,  \emph{The homology of a locally finite graph with ends}, Combinatorica 30 (2010), no. 6, 681--714
 

\bibitem{Georga} A. Georgakopoulos, \emph{Cycle decompositions: from graphs to continua}, Adv. Math. 229 (2012), no. 2, 935--967


\bibitem{Halmos} P.R. Halmos, \emph{Measure Theory}, The University Series in Higher Mathematics, New York: van Nostrand, 1950

\bibitem{SKHansen} S.K. Hansen, \emph{Measure homology}, Math. Scand., Vol. 83(1998), 205--219 

\bibitem{Hen} E. Henze, \emph{Einf\"uhrung in die Ma\ss theorie}, Mannheim:Bibliographisches Institut, 1971  

\bibitem{JGHocking} J. B. Hocking, G. S. Young, \emph{Topology}, London: Addison-Wesley Publishing Company inc., 1961

\bibitem{GMargulis} G. Margulis, \emph{Discrete subgroups of semisimple Lie groups}, Berlin: Springer Verlag, 1991

\bibitem{HJMunkholm} H. J., Munkholm, \emph{Simplices of maximal volume in hyperbolic space, Gromov's norm, and Gromov's proof of Mostow's rigidity theorem (following Thurston)}, Topology Symposium, Siegen 1979 (Proc. Sympos., Univ. Siegen, Siegen, 1979), pp. 109–124, Lecture Notes in Math., 788, Springer, Berlin, 1980

\bibitem{Preprint} J. Przewocki, \emph{Milnor-Thurston homology groups of the Warsaw Circle}, Topology Appl. 160 (2013), no. 13, 1732--1741

\bibitem{Ratcl} J.G. Ratcliffe, \emph{Foundations of Hyperbolic Manifolds},  Springer-Verlag: New York-Heidelberg-Berlin, 1994

\bibitem{HLRoyden} H. L. Royden, \emph{Real analysis}, 2nd ed., New York: Macmillan Publishing Co., Inc., London: Collier Macmillan Publishers, 1968

\bibitem{SrivastavaSM} S.M. Srivastava, \emph{A course on Borel Sets}, New York: Springer Verlag, 1998

\bibitem{Thurst} W.P. Thurston, \emph{Geometry and Topology of Three-manifolds}, Lecture notes, available at http://www.msri.org/publications/books/gt3m, Princeton, 1978

\bibitem{AZastrow} A. Zastrow, \emph{On the (non)-coincidence of Milnor-Thurston homology theory with singular homology theory}, Pacific J. Math. 186(1998) 369 -- 396

\end{thebibliography}
\end{document}